\documentclass{article} 

\usepackage{amsmath,amsfonts,bm}

\def\1{\bm{1}}

\def\vv{{\bm{v}}}

\def\vx{{\bm{x}}}
\def\vy{{\bm{y}}}

\DeclareMathAlphabet{\mathsfit}{\encodingdefault}{\sfdefault}{m}{sl}
\SetMathAlphabet{\mathsfit}{bold}{\encodingdefault}{\sfdefault}{bx}{n}

\def\gX{{\mathcal{X}}}

\def\der{{\mathrm{d}}}

\newcommand{\E}{\mathbb{E}}
\newcommand{\Ls}{\mathcal{L}}
\newcommand{\Is}{\mathcal{I}}
\newcommand{\R}{\mathbb{R}}

\DeclareMathOperator*{\argmin}{arg\,min}

\usepackage[left=1in,top=1in,right=1in,bottom=1in,letterpaper]{geometry}
\usepackage{hyperref}
\usepackage{url}
\usepackage{microtype}
\usepackage{graphicx}
\usepackage{subfigure}
\usepackage{booktabs} 
\usepackage{amsmath}
\usepackage{amsfonts}
\usepackage{amsthm}
\usepackage{xspace}
\usepackage{soul}
\usepackage{mathtools}
\usepackage{natbib}
\usepackage{breqn}
\usepackage[T1]{fontenc}
\usepackage[font=small,labelfont=bf,tableposition=top]{caption}
\usepackage[bottom]{footmisc}
\usepackage{tablefootnote}
\usepackage[dvipsnames]{xcolor}
\usepackage{float}

\DeclareCaptionLabelFormat{andtable}{#1~#2  \&  \tablename~\thetable}

\newtheorem{thm}{Theorem}

\newtheorem{remark}[thm]{Remark}

\newlength\mylength
\newlength\mylengthLong
\title{Joint Inference of Trajectory and Obstacle in Mean-Field Games via Bilevel Optimization}

\author{Han Huang \thanks{H. Huang (huangh14@rpi.edu) is with the department of mathematics, Rensselaer Polytechnic Institute. }\qquad 
Jiajia Yu \thanks{J. Yu (jiajia.yu@duke.edu) is with the department of mathematics, Duke University. } \qquad 
Tianyi Chen \thanks{T. Chen (tianyi.chen@cornell.edu) is with Department of Electrical and Computer Engineering,
Cornell University } \qquad 
Rongjie Lai \thanks{Corresponding author. R Lai (lairj@purdue.edu) is with the department of mathematics, Purdue University. }
}

\begin{document}
\date{}
\maketitle

\begin{abstract}

Mean field game (MFG) is an expressive modeling framework for systems with a continuum of interacting agents. While many approaches exist for solving the forward MFG, few have studied its \textit{inverse} problem. In this work, we seek to recover optimal agent trajectories and the unseen spatial obstacle given partial observation on the former. To this end, we use a special type of generative models, normalizing flow, to represent the trajectories and propose a novel formulation of inverse MFG as a bilevel optimization (BLO) problem. We demonstrate the effectiveness of our approach across various MFG scenarios, including those involving multi-modal and disjoint obstacles, highlighting its robustness with respect to obstacle complexity and dimensionality. Alternatively, our formulation can be interpreted as regularizing maximum likelihood trajectory learning with MFG assumptions, which improves generalization performance especially with scarce training data. Impressively, our method also recovers the hidden obstacle with high fidelity in this low-data regime.

\end{abstract}

\section{Introduction}
Mean-field games (MFG) models a wide variety of interactive systems with a large number of participating agents by prescribing an identical objective for all agents and seeking the optimal strategy over a time interval, a process analogous to solving for the Nash equilibrium in a finite $N$-player game \cite{MFG_game_thry}. Owing to its expressivity, MFG is instrumental in advancing understandings in economic\cite{MFG_econ1, MFG_econ2}, biological \cite{MFG_bio_decision_making}, robotic \cite{MFG_swarm_control}, and financial systems \cite{MFG_fin1, MFG_fin2, MFG_fin3}. As a result, it is of natural interest to consider inverse problems arising from MFG systems. Broadly speaking, we aim to infer hidden parameters in the MFG given partial observations of its solutions. In this work, we focus on MFGs where agents wish to arrive at a terminal location with minimal distanced traveled while avoiding spatial obstacles, a natural setup for a variety of problems in drone motion planning and animal behaviors. The inverse MFG problem is thus to recover the unknown obstacle given optimal trajectories, which is of interest in e.g., topography surveying.

In this work, we formulate the inverse MFG as a bilevel optimization problem and solve it with the penalty method to jointly obtain the latent obstacle as well as optimal trajectories. To the best of our knowledge, existing methods for inverse MFG needs measurement on the agent density value, which may be difficult to obtain or require an intermediate estimation step. In contrast, our framework is grounded on a trajectory-based MFG formulation and allows obstacle recovery based on observed trajectories that are samples of the underlying agent density. 

Crucially, our framework shows impressive empirical sampling efficiency in obstacle inference, a highly desirable feature since trajectory data may be costly to obtain. With fewer than a hundred trajectories, our method can wager an excellent guess at the latent obstacle. In addition, when the obstacle is well-approximated, the proposed learning framework is more cognizant of and consistent with the dynamical laws underlying the observed data. As a result, when training data are approximately optimal for a MFG system, trajectory learning can be made tractable even with very limited demonstrations. In this scarce data regime, the proposed framework generalize much better than the maximum likelihood alternative, so the bilevel training can be interpreted as an effective regularization that reduces overfitting in trajectory modeling.

To illustrate the efficacy of our approach, we conduct extensive experiments in various dimensions with disjoint, multi-modal obstacles. Our method consistently achieves approximately $0.1$ relative $L_2$ error in obstacle recovery, and the visualized results are qualitatively faithful to the ground truths. In addition, we show the model's impressive performance in the scarce data regime, highlighting its robustness in obstacle inference and regularizing effect in trajectory learning. Lastly, we showcase our method's advantage over traditional approaches by applying it to high dimensional inverse MFG, a scenario that is intractable for methods involving spatial discretizations. The proposed framework shows promise in synthetic examples, which can pave way for e.g., inference in configuration spaces.

\subsection{Related Works}

\paragraph{Mean Field Games}


Variational mean field games (MFG) extend the dynamic formulation of optimal transport (OT)~\cite{OT_book}, where the final density match is not strictly enforced and additional transport costs beyond kinetic energy can be considered. While traditional methods for solving OT and variational MFG problems are well-established in lower dimensions ~\cite{achdou2010mean, benamou2014augmented,benamou2017variational,benamou2000computational,jacobs2019solving,papadakis2014optimal,yu2021fast,cuturi2013sinkhorn,yu2023computational}, recent advances in machine learning have enabled the solving of high-dimensional MFGs. The work~\cite{NN_MFP} approximates the value function of high-dimensional deterministic MFGs using deep neural networks, \cite{MFG_NF} solves the MFGs using normalizing flow based on trajectory formula,  ~\citep{Tong_APACNet} solve high-dimensional MFGs in a stochastic setting using the primal-dual formulation of MFGs, which is trained similarly to generative adversarial networks~\cite{GAN,WGAN}. 

The goal of inverse MFG is to recover certain parameters in the MFG system given partial observations from its solutions. 
This is a relatively new field, where initial works proposed to recover running costs and interaction kernels from boundary measurements \cite{inverse_MFG_boundary} as well as agent density and velocity \cite{inverse_MFG_metric}. 
\cite{yu2024invmfg_blo} focuses on potential MFG and formulates the inverse MFG as a bilevel optimization problem. Different from our paper, the aforementioned works use an Eulerian coordinate and the measurements are mesh-based, which prevents the implementation in high dimensions.
\cite{guo2024decoding,zhang2025learning} integrate the framework of Gaussian process to infer costs and environments in MFGs. While their measurements can be mesh-free, the inference process requires collocation points which are analogues to a mesh. 
In contrast, our framework is completely mesh-free. 
It is based on a trajectory-based formulation of MFG and requires only demonstrations of optimal agent paths. In addition, our framework leverages network parametrizations to avoid spatial discretizations, making it tractable to conduct inference in high dimensional settings.

\paragraph{Normalizing Flows}

Normalizing Flows (NF) are a family of invertible neural networks that find applications in density estimation, data generation, and variational inference~\cite{NICE, RealNVP, NSF, Flow++, MAF, IAF, NAF}. In machine learning, NFs are considered deep generative models, which include variational autoencoders~\cite{VAE} and generative adversarial networks~\cite{GAN} as other well-known examples. The key advantage of NF over the alternatives is the exact computation of the data likelihood, which is made possible by its invertibility.

NF models impose a certain structure on the flow mapping to enable fast evaluation of the log-determinant of the Jacobian. For a non-exhaustive list of examples, RealNVP and NICE combine affine coupling flows to build simple yet flexible mappings~\cite{RealNVP, NICE}, while NSF and Flow$++$ use splines and a mixture of logistic CDFs as more sophisticated coupling functions~\cite{NSF, Flow++}. MAF and IAF are autoregressive flows with affine coupling functions~\cite{MAF, IAF}, while NAF and UMNN parameterize the coupling function with another neural network~\cite{NAF, UMNN}.

\paragraph{Inverse Reinforcement Learning}

The central problem in inverse reinforcement learning (IRL) is to infer the reward function and/or the policy of an expert agent given observational data on its behaviors \cite{IRL_survey, IRL_survey_2}. Typically, the interaction between the expert and the environment is modeled as a Markov decision process (MDP), and the expert's action are sampled from the optimal policy that solves the MDP. IRL offers two main applications. First, we can train agents to learn from  demonstration to act in accordance with the expert's preferences. This is referred to as behavioral cloning \cite{IRL_slides}, and examples include helicopter control \cite{IRL_heli_control}, parking lot navigation \cite{IRL_parking}, and simulated driving \cite{IRL_driving}. Second, the learned reward function serves as a useful basis for imparting expert knowledge, which is shown to be more transferable than the policy viz. changes in the environment \cite{IRL_transferability}. 

The connection between MFG and reinforcement learning (RL) is well-known and explored in~\cite{MFG_RL, MFG_RL_2, MFG_RL_3}. One can think of MFG as a particular flavor of RL with special forms for the environment and reward function. Namely, the agent action and the state transition model are encoded in an ODE system, and the reward is a linear combination of three cost terms. To our knowledge, the only work that connects MFG and IRL is \cite{IRL_MFG}, which proposes an algorithm for solving non-cooperative MFG systems based on the RL formalism. In contrast, our approach leverages the special structures in the MFG costs to parametrize and solve the MFG system with NF, which naturally decomposes the ODE flow and allows for efficient evaluation of the terminal cost.


\section{Background}
In this section, we provide background information on mean-field games and discuss solvers for the forward problem using normalizing flow. 
\subsection{Mean-Field Games}

The study of MFGs considers classic $N$-player games at the limit of $N\to \infty$~\cite{MFG_varMFG,huang2006large,huang2007large}. Assuming homogeneity of the objectives and anonymity of the players, we set up an MFG with a continuum of non-cooperative rational agents distributed spatially in $\mathbb{R}^d$ and temporally in $[0, T]$. For any fixed $t\in [0,T]$, we denote the population density of agents by $p(\cdot, t)$. For an agent starting at $\vx_0\in \R^d$, their position over time follows a trajectory $\vx: [0,T] \to \R^d$ governed by
\begin{equation} \label{agent_traj}
\left\{\begin{aligned}
    \der \vx(t) &= \vv(\vx(t),t)\der t, \quad \forall t\in [0,T]\\
    \vx(0) &= \vx_0,
\end{aligned}\right.
\end{equation}
where $\vv: \R^d \times [0,T] \to \R^d$ specifies an agent's action at a given time. For simplicity, we turn off the stochastic terms in~\eqref{agent_traj} in this work. Then, each agent's trajectory is completely determined by $\vv(\vx,t)$. To play the game over an interval $[t, T]$, each agent seeks to minimize their objective:
\begin{equation} \label{MFG}
\begin{aligned}
    J_{\vx_0,t }(\vv,p) &\coloneqq \int_t^T [L(\vx(s),\vv(\vx(s),s)) + I(\vx(s), p(\cdot, s)) ] \der s + M(\vx(T), p(\cdot, T))\\
    &\text{s.t. } \eqref{agent_traj} \text{ holds}.
\end{aligned}
\end{equation}
Each term in the objective denotes a particular type of cost. The running cost $L: \R^d \times \R^d \to \mathbb{R}$ is incurred by each agent's own action. A commonly used example for $L$ is the kinetic energy $L(\vx,\vv) = \|\vv\|^2$, which accounts for the total amount of movement along the trajectory. The running cost $I: \R^d \times \mathcal{P}(\R^d) \to \mathbb{R}$ is accumulated through each agent interacting with one another or with the environment. For example, this term can be an entropy that discourages the agents from grouping together, or a penalty for colliding with an obstacle~\cite{NN_MFP}. The terminal cost $M: \R^d \times \mathcal{P}(\R^d) \to \mathbb{R}$ is computed from the agents' final state, which typically measures a discrepancy between the final density $p(\cdot, T)$ and a desirable density $p_1 \in \mathcal{P}(\R^d)$.


The setup~\eqref{agent_traj}--\eqref{MFG} concerns the strategy of an individual agent. Under suitable assumptions, the work~\cite{MFG_varMFG} developed a macroscopic formulation of the MFG that models the collective strategy of all agents. Suppose there exist functionals $\mathcal{I}, \mathcal{M}: \mathcal{P}(\R^d) \to \mathbb{R}$ such that
\begin{align*}
    I(\vx,p) = \frac{\delta \mathcal{I}}{\delta p}(\vx), \quad M(\vx,p) = \frac{\delta \mathcal{M}}{\delta p}(\vx),
\end{align*}
where $\displaystyle\frac{\delta}{\delta p}$ is the variational derivative. Then, the functions $p(\vx,t)$ and $\vv(\vx,t)$ satisfying~\eqref{MFG} coincide with the optimizers of the following variational problem: 
\begin{equation}\label{var_MFG}
\begin{split}
    \inf_{p,\vv} J(p,\vv) \coloneqq \int_0^T \int_{\R^d} L(\vx, \vv(\vx,t),p(\vx,t) \der\vx \der t &+ \int_0^T \mathcal{I} (p(\cdot, t))\der t + \mathcal{M}(p(\cdot, T))\\
    s.t. \quad \partial_t p(\vx,t) + \nabla_{\vx} \cdot (p(\vx,t)\vv(\vx,t)) &= 0, \quad \vx\in \R^d, t\in [0,T]\\
    p(\vx,0) = p_0(\vx), \quad \vx\in &\R^d.
\end{split}
\end{equation}
This formulation is termed the \emph{variational MFG}. Unless specified otherwise, we assume $T=1$ and use the $L_2$ transport cost $ L(\vx,\vv(\vx,t),p(\vx,t)) = \lambda_L p(\vx,t) \|\vv(\vx,t)\|_2^2$ from now on, where $\lambda_L \geq 0$ is as a hyperparameter. In addition, we consider interaction costs of the form $\mathcal{I} (p(\cdot, t)) = \int_{\R^d} B(\vx)p(\vx,t)dx$, where $B: \mathbb{R}^d \to \mathbb{R}$ may be thought of as an obstacle that penalizes agents who travel in its proximity.

\subsection{Trajectory-Based MFG}
An equivalent formulation of \eqref{var_MFG} is proposed in~\cite{MFG_NF} and is more suitable for our purpose as it accentuates the role of agent trajectories. Let $P(\cdot,t)$ be the measure that admits $p(\cdot,t)$ as its density for all $t\in [0,T]$. Define the agent trajectory as $F: \mathbb{R}^{d} \times \mathbb{R} \to \mathbb{R}^d$, where $F(\vx,t)$ is the position of the agent starting at $\vx$ and having traveled for time $t$. It satisfies the following differential equation:
\begin{equation} \label{NF_traj}
\left\{\begin{aligned}
    \partial_t F(\vx,t) &=  \vv(F(\vx,t), t), \quad \vx\in \R^d, t\in[0,T] \\
    F(\vx,0) &= \vx, \quad \vx\in \R^d.
\end{aligned}\right.
\end{equation}
The evolution of the population density is determined by the movement of agents. Thus, $P(\cdot,t)$ is simply the push-forward of $P_0$ under $F(\cdot, t)$; namely, $P(\cdot, t) = F(\cdot, t)_* P_0$, whose associated density satisfies
\begin{equation} \label{density_push-forward}
\begin{aligned}
    p(\vx,t) &= \mathrm{d}(F(\cdot, t)_* P_0)(\vx),
\end{aligned}
\end{equation}
where $\mathrm{d}(F(\cdot, t)_* P_0(\vx))$ is a Radon-Nikodym derivative. Hence, \eqref{var_MFG} becomes:

\begin{equation} \label{var_MFG_traj}
\begin{aligned}
    \inf_{F} \quad \lambda_L \int_0^1 \int_{\R^d}  \|\partial_t F(\vx,t)\|_2^2 p_0(\vx)\der\vx \der t &+ \lambda_{\mathcal{I}} \int_0^1 \int_{\R^d} B(F(x,t))p_0(\vx)\der\vx \der t + \lambda_{\mathcal{D}}  \mathcal{D} (P_1, F(\cdot, T)_*P_0) \coloneqq \mathcal{L} (F; B)\\
    s.t. \quad  &F(\vx,0) = \vx.
\end{aligned}
\end{equation}
Where $D$ is a discrepancy between two distributions.

\subsection{Normalizing Flows}


Suppose we have a dataset $\gX= \{\vx_n\}_{n=1}^N \subseteq \mathbb{R}^d$ generated by an underlying data distribution $P_1$; that is, $\vx_n \overset{\mathrm{iid}}{\sim} P_1$. We define a \textit{flow} to be an invertible function $f_\theta:\mathbb{R}^d \to \mathbb{R}^d$ parameterized by $\theta \in \mathbb{R}^W$, and a \textit{normalizing flow} to be the composition of a sequence of flows: $F_\theta = f_{\theta_K} \circ f_{\theta_{K-1}} \circ ... \circ f_{\theta_2} \circ f_{\theta_1}$ parameterized by $\theta = (\theta_1, \theta_2, ..., \theta_K)$. To model the complex data distribution $P_1$, the idea is to use an NF to gradually transform a simple base distribution $P_0$ to $P_1$. Formally, the transformation of the base distribution is conducted through the push-forward operation $F_{\theta*}P_0$, where  $F_{\theta*}P(A) \coloneqq P(F_{\theta}^{-1}(A))$ for all measurable sets $A\subset\R^d$. The aim is that the transformed distribution resembles the data distribution. Thus, a commonly used loss function for training the NF is the KL divergence:
\begin{equation}\label{KL_NF}
    \min_\theta \quad  \mathcal{D}_{KL} (P_1 || F_{\theta*}P_0).
\end{equation}
If the measure $P$ admits density $p$ with respect to the Lebesgue measure, the density $F_{\theta*}p$ associated with the push-forward measure $F_{\theta*}P$ satisfies the change-of-variable formula~\cite{NF_survey}:
\begin{equation}\label{change_of_var}
    F_{\theta*}p(\vx) = p(F_\theta^{-1}(\vx))|\det \nabla F_\theta^{-1}(\vx)|,
\end{equation}
where $\nabla F_\theta^{-1}$ is the Jacobian of $F_\theta^{-1}$. This allows one to tractably compute the KL term in~\eqref{KL_NF}, which is equivalent to the negative log-likelihood of the data up to a constant (independent of $\theta$):
\begin{equation} 
\begin{split}
    \mathcal{D}_{KL} (P_1 || F_{\theta*}P_0) &= - \mathbb{E}_{\vx\sim P_1} [\log F_{\theta*}p_0(\vx)] + \text{const.} \\
    &= -\mathbb{E}_{\vx\sim P_1} \left[\log p_0(g_{\theta_1} \circ ... \circ ... g_{\theta_K}(\vx)) + \sum_{k=1}^K \log |\det \nabla g_{\theta_k}(\vy_{K-k})|\right] + \text{const.},
\end{split}
\label{NF_NLL}
\end{equation}
where $g_{\theta_i} \coloneqq f_{\theta_i}^{-1}$ are the flow inverses, and $\vy_j \coloneqq g_{\theta_{K-j}} \circ ... \circ g_{\theta_K}(\vx), \forall j = 1, 2, ..., K-1, \vy_0 \coloneqq \vx$.

We remark that there exists flexibility in choosing the base distribution $P_0$, as long as it admits a known density to evaluate on. In practice, a typical choice is the standard normal $\mathrm{N}(0,I)$. In addition, \eqref{NF_NLL} requires the log-determinant of the Jacobians, which take $O(d^3)$ time to compute in general. Existing NF architectures sidestep this issue by using parametrizations that follow a (block) triangular structure~\cite{NF_survey}, so that the log-determinants can be computed in $O(d)$ time from the (block) diagonal elements.


\section{Methodology}

\subsection{Problem Setup}

In this work, we assume a group of agents initially distributed according to $P_0$. Each agent act optimally according to its preferences specified by the MFG system with an unseen obstacle $B(\vx)$, giving rise to a trajectory $\vx(t)$ satisfying \eqref{MFG} with $\vx(0) = \vx_0$ as the initial position. In addition, we assuming the underlying MFG system admits a variational form \eqref{var_MFG}, and agents traverse the shortest path to arrive at their destination specified as $P_1 \in \mathcal{P}(\R^d)$ with terminal fidelity given by Jeffery's divergence weighted by $\lambda_{\mathcal{D}} >0$: $\mathcal{M}(P(\cdot, T)) = \lambda_{\mathcal{D}} \mathcal{D} (P(\cdot, T), P_1)= \lambda_{\mathcal{D}} [\mathcal{D}_{KL}(P(\cdot, T) || P_1) + \mathcal{D}_{KL}(P_1|| P(\cdot, T) )]$. This is a natural setup for e.g., motion planning, where agents can be drones  and trajectory values are spatial coordinates. In this case, solving the MFG with $L_2$ transport cost imparts the preference for drones to arrive at their designated locations following the shortest path, and at the same time avoid any obstacles in the interim.

 Define $P_D$ as the distribution of optimal agent trajectories according to the MFG system \eqref{MFG}, i.e., $\vx(t) \sim P_D$. As training data, we have access to $N$ sample trajectories represented by their values at $K$ regularly spaced temporal junctures: $D = \{\vx^i\}_{i=1}^N \subset \mathbb{R}^{K \times d}, \vx^i_k \coloneqq \vx^i(t_k), t_k = k\Delta t, \Delta t = 1/K$. Our goal here is twofold: first, we want to solve the MFG system, in the sense that given a new agent starting at $\hat{\vx}_0$, we can deduce (values on) its optimal trajectory $\hat{\vx}^i(t_k) \coloneqq \hat{\vx}^i_k$ satisfying \eqref{MFG}. In the case of drone motion planning, this translates to being able to send new drones on courses that resemble training examples, i.e., following the shortest path while sidestepping obstacles.
 
 Second, we want to infer the obstacle $B(\vx)$ present in the MFG system, which can be considered as an inverse problem where the agents exploring the surroundings play the role of the probe to uncover the unknown environment. In drone-based applications, this is central to the task of automated surveying \cite{drone_survey_buildings}. It turns out the two tasks can be achieved simultaneously by solving a bilevel problem that couples unoptimized agent trajectories with unseen obstacles. Before that, however, we first need to reformulate the variational problem \eqref{var_MFG} to one that is more explicit in terms of agent trajectories.

\subsection{Bilevel Formulation of an Inverse MFG Problem}

We now derive the central optimization problem this work seeks to solve. We know that any given obstacle $B$ specifies a forward MFG system \eqref{var_MFG_traj} and results in optimal agent trajectories $F^*$. Among all possible obstacles, we look for one that produces  trajectories closest to the training examples, yielding the following bilevel optimization problem:
\begin{equation} \label{BLO}
\begin{aligned}
    \min_{B, F} &\quad \E_{\vx(t) \sim P_D, t \sim \text{Unif}[0,1]} d(F^*(\vx(0),t; B), \vx(t)) \\
    s.t. &\quad F^* = \argmin_{F(\vx,0) = \vx} \Ls(F; B),
\end{aligned}
\end{equation}
where $F^*(\vx(0),t; B)$ denotes the optimal trajectory for a agent starting at $\vx(0)$ for the MFG problem with obstacle $B$, and $d(\cdot, \cdot)$ is the criterion for determining the closeness between two trajectories. From here on, we consider $d(\vx,\vy) = \|\vx-\vy\|^2_2$. As we have access to $K$ equidistant temporal values of samples $\vx(t) \sim P_D$, we consider the empirical rule with temporal discretization $t_k = k \Delta t, \Delta t = K$, resulting in the following problem:
\begin{equation} \label{BLO_ERM}
\begin{aligned}
    \min_{B, F} &\quad \frac{1}{KN} \sum_{n=1}^{N} \sum_{k=1}^{K} \| F^*_k(\vx^n_0; B) - \vx^n_k\|^2_2 \\
    s.t. &\quad F^* = \argmin_{F(\vx,0) = \vx} \Ls(F; B),
\end{aligned}
\end{equation}
where $F^*_k(\vx; B) \coloneqq F^*(\vx, t_k; B)$.

\subsection{Parametrization of Agent Trajectories and Obstacles}
As it stands, the empirical rule \eqref{BLO_ERM} presents an infinite-dimensional problem that cannot be tractably solved. 
We first conduct a time discretization of the trajectory formula \eqref{var_MFG_traj} for the lower level MFG loss function. 
Since the terminal cost is invariant to discretization, we derive the discrete versions for the transport and interaction costs in \eqref{var_MFG_traj}. Let $ t_k \coloneqq k \cdot \Delta t, \forall k = 0,1,...K, K=\frac{1}{\Delta t}$ define a regular grid in time. Denote $F(\vx,t_k) \coloneqq F_k(\vx)$. The transport cost becomes:
\begin{equation}
\begin{aligned}
 \lambda_L \int_0^1 \int_{\R^d}  \|\partial_t F(\vx,t)\|_2^2 p_0(\vx)\der\vx \der t &= \lambda_L \E_{x\sim P_0} \int_0^1  \|\partial_t F(\vx,t)\|_2^2 \der t \\
 &\sim \lambda_L \E_{x\sim P_0} \sum_{k=0}^{K-1} \|\frac{F(\vx,t_{k+1}) - F(\vx,t_{k})}{\Delta t}\|_2^2 \Delta t\\
 &\sim \frac{\lambda_L K}{M} \sum_{m=1}^M \sum_{k=0}^{K-1} \|F_{k+1}(x^m) - F_{k}(x^m)\|_2^2 
\end{aligned}
\end{equation}
where we used forward Euler discretization on the second line, then applied the empirical rule at the last step. Parametrizing each $F_k$ as $F_{\theta_k}$ gives the objective in \eqref{var_MFG_disc}.
A similar procedure is used for deriving the discretized interaction cost:
\begin{equation}
\begin{aligned}
 \lambda_{\mathcal{I}} \int_0^1 \int_{\R^d} B(F(x,t))p_0(\vx)\der\vx \der t &=  \lambda_{\mathcal{I}} \E_{x\sim P_0} \int_0^1  B(F(x,t)) \der t \\
 &\sim \frac{\lambda_{\mathcal{I}}}{K} \E_{x\sim P_0} \sum_{k=0}^{K-1} B(F_{k}(x))\\
 &\sim \frac{\lambda_{\mathcal{I}}}{KM} \sum_{m=1}^M \sum_{k=0}^{K-1} B(F_{k}(x^m))
\end{aligned}
\end{equation}

Following \cite{MFG_NF}, we further introduce appropriate parametrizations for both agent trajectories as well as the obstacle. Define flow maps for $0\leq t_i \leq t_j \leq 1$ as  $\Phi_{t_{i}}^{t_{j}} (\vx_0) = \vx(t_{j})$, where
\begin{equation}\label{flow_map}
\left\{\begin{aligned}
    \der\vx(t) &= \vv(\vx(t),t)\der t, \quad t\in [t_i, t_j]\\
    \vx(t_i) &= \vx_0.
\end{aligned}\right.
\end{equation}
 we leverage the semi-group property $\Phi_{t_{b}}^{t_{c}} \circ \Phi_{t_{a}}^{t_{b}} = \Phi_{t_{a}}^{t_{c}}$ to decompose agent trajectories into a sequence of flow maps.  It follows that $F_k = \Phi_{t_{k-1}}^{t_k} \circ \Phi_{t_{k-2}}^{t_{k-1}} \circ ... \circ \Phi_{t_0}^{t_1}, \forall k = 1, 2, ... K$, and we parameterize each $ \Phi_{t_i-1}^{t_i}: \mathbb{R}^d \to \mathbb{R}^d$ by a normalizing flow architecture $f_{\theta_i}$ (detailed flow structure is provided in numerical part) . In addition, we parametrize the obstacle with a simple MLP: $B(\vx) = B_\phi(\vx)$. Denote $F_{\theta_k} \coloneqq f_{\theta_k} \circ ... \circ f_{\theta_1}$, $F_{\theta_0} \coloneqq Id$, and $\theta = (\theta_1, ..., \theta_K)$. 
Therefore, the flow-parametrized bilevel problem becomes:
\begin{equation} \label{BLO_ERM_parametrized}
\begin{aligned}
    \min_{\phi, \theta} &\quad \frac{1}{KN} \sum_{n=1}^{N} \sum_{k=1}^{K} \| F_{\theta_k^*}(\vx^n_0) - \vx^n_k\|^2_2\\
    s.t. &\quad \theta^* = \argmin_{\hat{\theta}} \hat{\Ls}(\hat{\theta}; \phi).
\end{aligned}
\end{equation}
where 
\begin{equation}\label{var_MFG_disc}
\begin{aligned}
\hat{\Ls}(\hat{\theta}; \phi) = \frac{\lambda_L K}{M}\cdot  \sum_{m=1}^M \sum_{k=0}^{K-1}\|F_{\hat{\theta}_{k+1}}(\vx^m) - F_{\hat{\theta}_k}(\vx^m)\|_2^2 + \frac{\lambda_\Is}{KM} \sum_{m=1}^M \sum_{k=0}^{K-1} B_\phi(F_{\hat{\theta}_k}(\vx^m)) + \lambda_{\mathcal{D}} D (P_1, F_{\hat{\theta}_K *}P_0)
\end{aligned}
\end{equation}
is the discretized empirical rule of \eqref{var_MFG_traj} for $\vx^m \sim P_0$. Note that the initial identity flow is automatically satisfied, so the lower level problem is unconstrained. 

\begin{remark}
    The upper objective function in ~\eqref{BLO_ERM_parametrized} can be equivalently derived from the maximum likelihood principle by assuming Markovian trajectories and using Gaussian transition kernels. Suppose we have trajectory data $\{x^n\}^N_{n=1} \sim P_D$, $x^n = (x^n_1, ..., x^n_K) \in \mathbb{R}^{d \times K}$. The maximum likelihood principle seeks a model parametrized by $\theta$ that is most likely to produce the observed trajectories:

\begin{equation} \label{BLO_MLE_obj}
\begin{aligned}
    \max_\theta \quad \mathbb{E}_{x \sim P_D} [\log p_\theta(x)].
\end{aligned}
\end{equation}

To make this optimization tractable, we assume the trajectories are Markovian, i.e., $p_\theta(x_k | x_1, x_2, ..., x_{k-1}) = p_\theta(x_k | x_{k-1}), k=1, ..., K$, and decompose the likelihood as: $p_\theta(x) = p_0(x_0) \prod_{k=0}^{K-1}p_{\theta_{i+1}}(x_{k+1}|x_k)$. We then parametrize each transition density as a conditional gaussian: $p_{\theta_{k+1}}(x_{k+1}|x_k) = \mathrm{N}(x_{k+1}; f_{\theta_{k+1}}(x_k), \sigma^2I)$, where each $f_{\theta_{k+1}}$ is a flow network. With these in mind, the maximum likelihood objective~\eqref{BLO_MLE_obj} becomes:

\begin{equation} \label{BLO_MLE_flow_obj}
\begin{aligned}
    &\max_\theta \quad \mathbb{E}_{x \sim P_D} [\log p_\theta(x)]\\
    = &\max_\theta \quad \mathbb{E}_{x \sim P_D} [\log p_0(x_0) + \sum_{k=0}^{K-1}\log p_{\theta_{k+1}}(x_{k+1}|x_k)]\\
    \sim &\max_\theta \quad \sum_{k=0}^{K-1} \mathbb{E}_{x \sim P_D} [\|f_{\theta_{k+1}}(x_k) - x_{k+1}\|^2_2]\\
    = &\max_\theta \quad \sum_{k=1}^K \E_{x \sim P_D} [\|F_{\theta_k}(x_0) - x_k\|^2_2]. 
\end{aligned}
\end{equation}
Here, applying the empirical rule on the expectation yields the same objective as ~\eqref{BLO_ERM_parametrized}.

\end{remark}

\subsection{Penalty-based Algorithm for Bilevel Optimization}

To solve \eqref{BLO_ERM_parametrized}, we first reformulate the bilevel problem to an equivalent single-level problem with an inequality constraint. 

\begin{equation} \label{BLO_single_with_constraint}
\begin{aligned}
    \min_{\phi, \theta} &\quad \frac{1}{KN} \sum_{n=1}^{N} \sum_{k=1}^{K} \| F_{\theta_k}(\vx^n_0) - \vx^n_k\|^2_2 \\
    s.t. &\quad H(\phi) \geq \hat{\Ls}(\theta; \phi),
\end{aligned}
\end{equation}

where $H(\phi) \coloneqq \min_\theta \hat{\Ls}(\theta; \phi)$. Observe that any feasible $(\theta, \phi)$ in \eqref{BLO_single_with_constraint} must satisfy $\theta = \argmin_{\hat{\theta}} \hat{\Ls}(\hat{\theta}; \phi)$, so its optimizers coincide with the original problem \eqref{BLO_ERM_parametrized}. We use the penalty method to solve this constrained problem:
\begin{equation} \label{BLO_penalty}
\begin{aligned}
    \min_{\phi, \theta} &\quad \frac{1}{KN} \sum_{n=1}^{N} \sum_{k=1}^{K} \| F_{\theta_k}(\vx^n_0) - \vx^n_k\|^2_2 + \lambda_P [\hat{\Ls}(\theta; \phi) - H(\phi)]_+, \\
\end{aligned}
\end{equation}


where $\lambda_p >0$ is a large penalty weight picked as a hyperparameter. Assuming $H, \hat{\Ls}$ are continuously differentiable, we can leverage the envelop theorem to compute the gradient $\nabla_\phi H(\phi) = \nabla_\phi \hat{\Ls}(\theta^*;\phi)$ for any $\theta^*\in \argmin_{\hat{\theta}} \hat{\Ls}(\hat{\theta}; \phi)$ directly via auto-differentiation \cite{envelop_thm}. The problem \eqref{BLO_penalty} can thus be solved with any black-box optimization algorithm. Crucially, even though $H(\phi)$ is obtained via an optimization procedure, our approach not need to backpropagate through the scheme to obtain hypergradients. As such, we  sidestep the linear growth in memory consumption viz. the length of optimization dynamics. 

To determine whether the penalty method is the ideal algorithm for problem~\eqref{BLO_ERM_parametrized}, we  explored alternative approaches for bilevel optimization such as the gradient method \cite{BLO_survey}, BDA \cite{BDA}, IAPTT-GM \cite{IAPTT}, and found the penalty-based framework performing favorably in both speed and robustness. We also experimented with the augmented lagrangian method, but it is more expensive as the function $H(\phi)$ had to be approximated at each inner step of the primal sub-problem. 

In addition, we remark that equation~\eqref{BLO_penalty} has the form of a maximum a posteriori (MAP) estimation, where the trajectory likelihood term is regularized by a prior that reflects our assumption that the observed data are solutions of a MFG system. Empirically, as we will see later, this regularization is instrumental for reducing overfitting in trajectory inference, especially when the amount of training data available is limited.

\subsection{Identifiability and Regularization}

In statistical inference, the concept of identifiability \cite{identifiability} is concerned with whether the quantities of interest can be uniquely determined from given data. For our setup, we distinguish between weight space and physical space identifiability.
\paragraph{Weight Space}
Weight space identifiability amounts to whether \eqref{BLO_ERM_parametrized} admits a unique solution $(\theta^*, \phi^*)$. In general, the answer is negative given the overparametrized nature of neural networks. For example, consider a MLP with consecutive layer weights and biases $(W_i, b_i), (W_{i+1}, b_{i+1})$. An equivalent MLP can be obtained by applying a permutation to the rows of $W_i$ and entries of $b_i$, then applying the same permutation to the columns of $W_{i+1}$ and entries of $b_{i+1}$ \cite{BNN}. However, in practice we are not concerned with identifying the network weights as they carry no physical meaning. Instead, our desiderata is for the networks to represent the underlying agent trajectories and obstacles accurately. Additionally, this symmetry poses no issue for our optimization, perhaps unsurprising given that MLPs can be trained to convergence under typical settings.

\paragraph{Physical Space}
Physcial space identifiability is defined as the uniqueness of the optimizers $(B^*, F^*)$ for \eqref{BLO_ERM}, which are the actual entities of interest. If we are given any solution $(B^*, F^*)$ for \eqref{BLO_ERM}, another solution may be obtained by simply offsetting $B^*$ with a constant: $(\Tilde{B}^*, F^*) = (B^*+c, F^*)$, which maintains validity on the lower-level. To mitigate this non-uniqueness, we follow typical inverse problem approaches to introduce a regularization. The idea is to match the total mass between the true and parametrized obstacles:
\begin{equation} \label{mass_reg}
\begin{aligned}
    \mathcal{R}(B) = [\frac{\int_\Omega B(\vx)\der \vx - \int_\Omega B^*(\vx)\der \vx}{\int_\Omega B^*(\vx)\der \vx}]^2
\end{aligned}
\end{equation}

We introduce parametrization into the regularization by replacing $B(\vx) \to B_\phi(\vx)$ to get $\mathcal{R}(\phi)$. The final optimization problem is

\begin{equation} \label{BLO_penalty_mass_reg}
\begin{aligned}
    \min_{\phi, \theta} &\quad \frac{1}{KN} \sum_{n=1}^{N} \sum_{k=1}^{K} \| F_{\theta_k}(\vx^n_0) - \vx^n_k\|^2_2 + \lambda_P [\hat{\Ls}(\theta; \phi) - H(\phi)]_+ + \lambda_M \mathcal{R}(\phi), \\
\end{aligned}
\end{equation}

We can do this without loss of generality, because the ground truth obstacle mass can always be normalized to, or regarded as, one. Empirically, mass regularization helps stabilize our optimization procedure.

\section{Numerical Results}

We empirically demonstrate the validity of our approach by formulating and solving the bilevel problem \eqref{BLO_ERM_parametrized}. Our experiments underscore the exceptional robustness of our method as we apply it across a range of MFG  scenarios, encompassing disjoint and multi-modal obstacles in varying dimensions. Moreover, our framework consistently delivers strong obstacle recovery results while significantly improving the efficiency of dynamics learning, particularly in situations where data is scarce.

For all synthetic examples, we first solve the discretized forward problem \eqref{var_MFG_traj} with flow parametrizations, then use it to sample agent trajectories as training data for the bilevel problem \eqref{BLO_ERM_parametrized}. Unless specified otherwise, all forward problems are optimized until convergence with the NSF \cite{NSF} coupling flow and Adam \cite{ADAM}, a setting that is widely applicable with minimal tuning. 


\begin{figure}[th]
\centering
\begin{minipage}{0.49\linewidth}
\includegraphics[width=.85\linewidth]{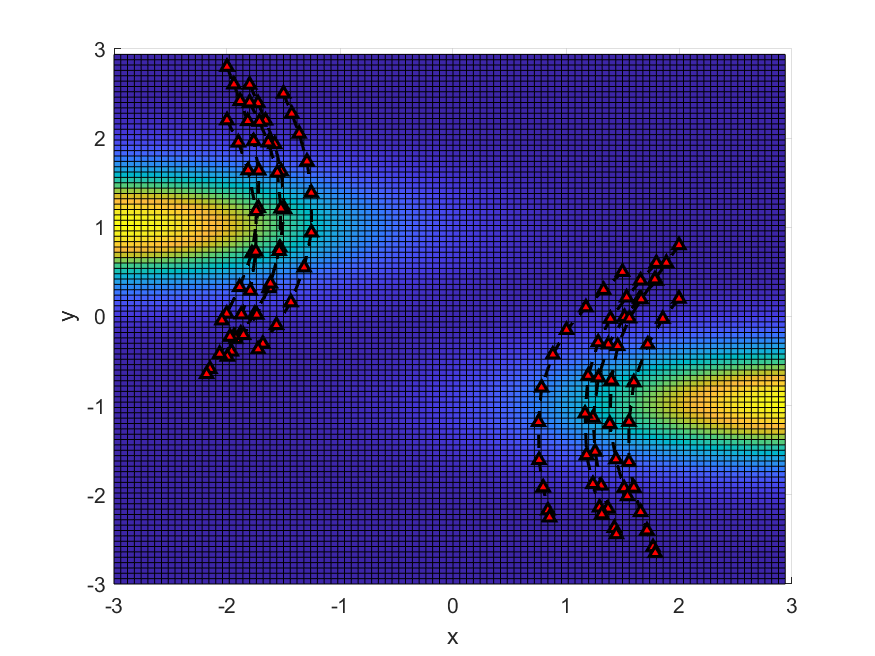}\\
\end{minipage}\hfill
\centering
\begin{minipage}{0.49\linewidth}
\includegraphics[width=.85\linewidth]{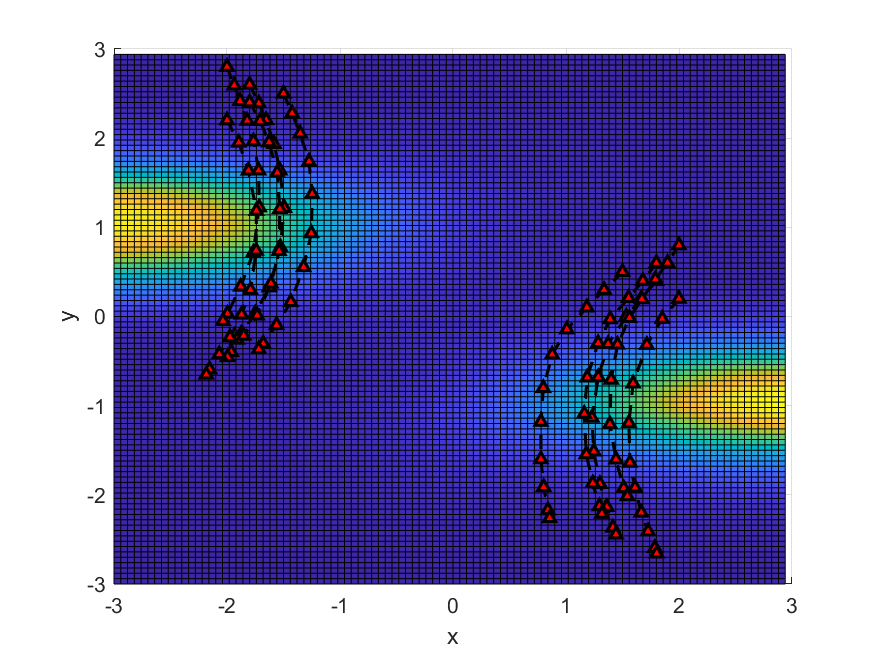}\\
\end{minipage}\hfill
\vskip -10pt
\centering
\begin{minipage}{0.49\linewidth}
\includegraphics[width=.85\linewidth]{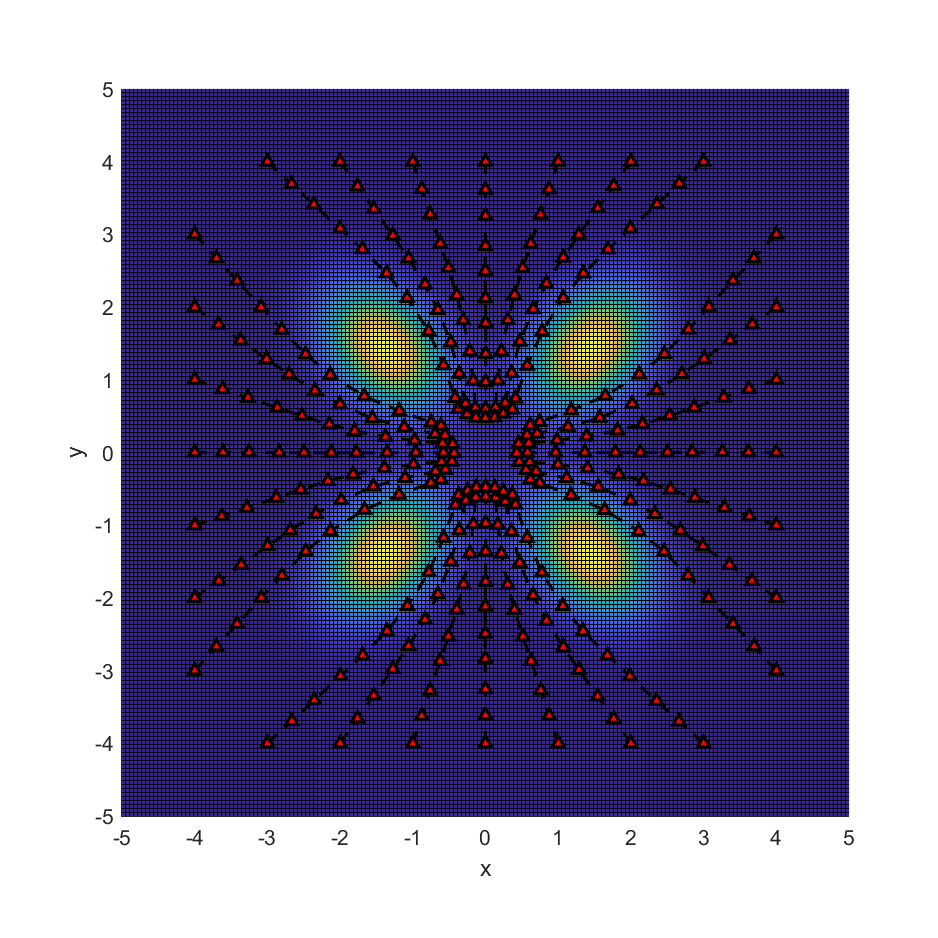}\\
\end{minipage}\hfill
\centering
\begin{minipage}{0.49\linewidth}
\includegraphics[width=.85\linewidth]{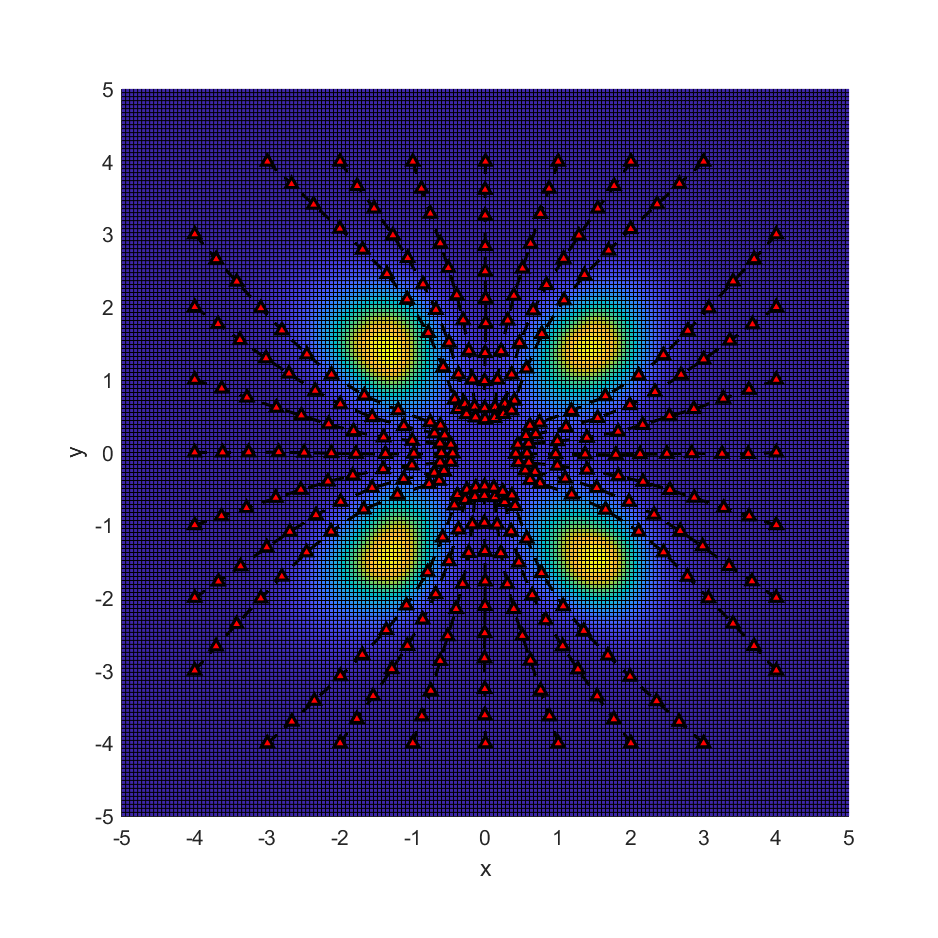}\\
\end{minipage}\hfill
\vskip -10pt
\caption{Left: Ground truth two-bars(top) and flower(bottom) obstacle with sample trajectories overlaid; Right: inferred obstacle from trajectories with recovered sample trajectories overlaid.}
\label{Two_bars_BLO_results}
\end{figure}

\paragraph{Detailed Experimental Setup}
\label{supp:numerical_experiments}

We use a popular normalizing flow architecture, NSF-CL \cite{NSF}, to parametrize agent trajectories. Each NSF-CL flow step consists of an alternating block of rational-quadratic splines as the coupling function followed by a linear transformation. The conditioner has 256 hidden features, 8 bins, 2 transform blocks, ReLU activation, and a dropout probability of 0.25. 

We parametrize the obstacles as a simple fully connected network. It has 3 layers with width 128 as well as residual connections. ReLU activation is used, though we have also experimented with Softmax, tanh, elu \cite{Elu}, and Mish \cite{Mish}. It is worth remarking that ReLU should not be used for activation if one needs to compute the hypergradient, because it is zero almost everywhere after differentiating twice. This behavior has no bearing on the penalty method, another reason for why it is preferred. An exponential activation is used after the last layer to ensure positivity in the network's output, which we empirically observe to expedite the optimization. 

In all but the scarce data experiments, we train the forward problem with Adam \cite{ADAM} at a $10^{-3}$ learning rate until convergence, then sample $3\cdot 10^4$ agent trajectories as training data and $10^4$ as testing data. A pretrained model on the same dataset with the MLE objective \eqref{BLO_MLE_obj} for $2\cdot 10^4$ steps is used as initialization for the BLO problem, which is trained for $10^4$ optimization steps via Adam at learning rates of $3\cdot 10^{-4}, 10^{-2}$ for trajectories and obstacles, respectively. In addition, a scheduler decays the learning rates by a factor of 0.8 every $10^3$ steps. In each outer optimization step, the approximation of lower level optimality in \eqref{BLO_penalty} is done via 50 gradient steps, also through Adam at a step size of $10^{-5}$. 

For the mass regularization, we use a regular grid in space and the Simpson rule to approximate the mass integral with 30 grid points in each dimension. This proves sufficient in two and three dimensional settings. For higher dimensional setups, one can opt for a Monte-Carlo or Quasi Monte-Carlo method\cite{Quasi_monte_carlo}, where the latter helps reduce the estimation variance.

For convenience, we summarize the relative $L_2$ obstacle errors for all experiments in table~\ref{table:obs_rel_l2_err} and present the qualitative results in the following subsections. Our approach consistently achieves $\sim0.1$ relative error with minimal tuning on training hyperparameters, e.g., $\lambda_P, \lambda_M,$ and learning rates, demonstrating robustness in the presence of disjoint, multi-modal obstacles and varying dimensionality. 

\begin{table}[h]
\begin{center}
\begin{small}
\begin{sc}
\begin{tabular}{lccccccc}
\toprule
Dataset & Two Bars & Flower & Cylinders & Castle & Mountains \\
\midrule
Rel. $L_2$ err   & 0.047 & 0.084 & 0.148 & 0.064 & 0.082 &  \\
\bottomrule
\end{tabular}
\end{sc}
\end{small}
\end{center}
\caption{Relative $L_2$ errors of recovered obstacles for different MFG settings.}
\label{table:obs_rel_l2_err}
\end{table}

\subsection{2D Examples}

We first study MFG defined in $\mathbb{R}^2$, where agents move on a two dimensional plane and seek to avoid obstacles in this space. We evaluate our method on two sets of obstacles with increasing complexity, demonstrating accurate recoveries of obstacle and trajectory in both cases.

\paragraph{Two Bars}

Shown in Figure \ref{Two_bars_BLO_results}, we consider an obstacle that has two connected components and modes: $B(\vx) = 50\cdot \frac{1}{2}\sum_{i=1}^2 \mathrm{N}(\vx; \mu_i, \Sigma_i), $ $\mu_1 = (-2, 1.25), \mu_2 = (4, 0.75), \Sigma_1 = \Sigma_2 = \text{diag}(1.5, 0.01)$, and the detectors are gaussian mixtures below and above the obstacles. Note that the domain of interest is picked to be $\Omega = [-3,3]^2$, so the obstacle has more support outside of $\Omega$. 

\paragraph{Flower}

In Figure \ref{Two_bars_BLO_results}, we present an example with more disjoint parts. The obstacle in question is a mixture of identical copies of the same gaussian in four quadrants: $\mathrm{N}(\vx; \mu, \Sigma), $ $\mu = (\frac{2\sqrt{2}}{2}, \frac{2\sqrt{2}}{2}), \Sigma = \text{diag}(0.2, 0.4)$. The agents start as a gaussian centered at the origin and travel radially outwards, providing enough coverage to recover the underlying obstacle. \\



In both scenarios, we overlay the obstacle and trajectories jointly obtained from solving \eqref{BLO_ERM_parametrized} with $N=30k$ training data and compare them to the ground truth obstacle and trajectories. Our approach is effective in accurately capturing the obstacle and modeling the underlying agent dynamics, which is evidenced by a close resemblance of obstacle avoidance behaviors observed in both the ground truth and inferred trajectories.

\subsection{3D Examples}

We then turn to MFG defined in $\mathbb{R}^3$, which is suitable for many motion planning problems carried out by animals and drones. Similar to before, we evaluate our method against a sequence of increasingly difficult setups to highlight its effectiveness.

\paragraph{Three Cylinders}

\begin{figure}[h]
\vskip -10pt
\centering
\begin{minipage}{0.5\linewidth}
\includegraphics[width=.8\linewidth]{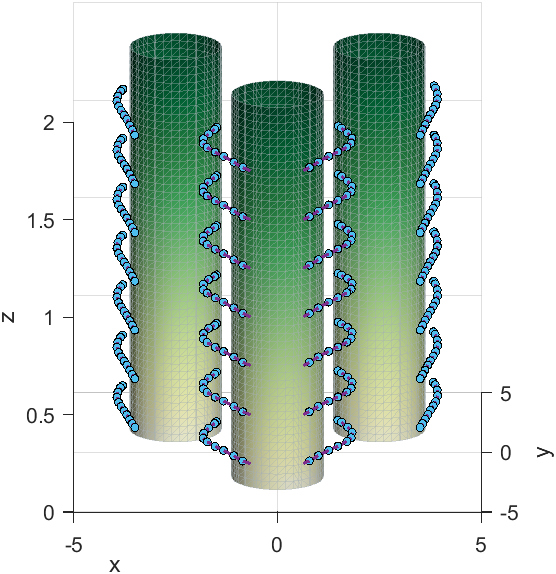}\\
\end{minipage}\hfill
\begin{minipage}{0.5\linewidth}
\includegraphics[width=.8\linewidth]{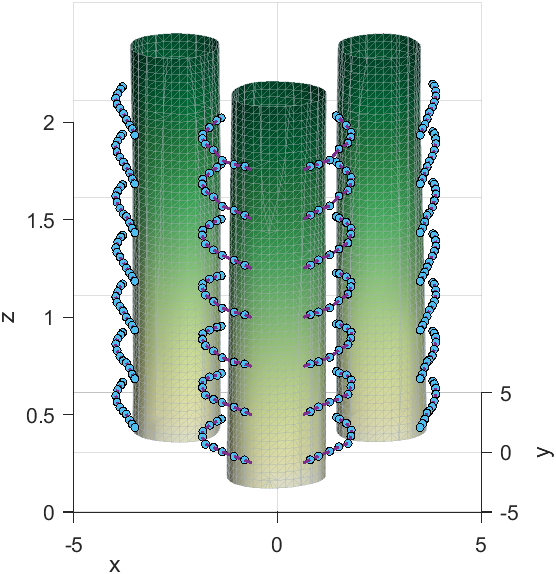}\\
\end{minipage}\hfill
\caption{Left: ground truth cylinder obstacle with sample trajectories overlaid; Right: inferred obstacle from trajectories with predicted trajectories. The obstacle color gradient indicates its height. }
\label{3Cyl_BLO_results}
\end{figure}

We start with a simple constant extension of 2D gaussian mixture: $50\cdot \frac{1}{3}\sum_{i=1}^3 \mathrm{N}(\vx; \mu_i, \Sigma_i), $ $\mu_1 = (0, -2), \mu_2 = (2.5, 2), \mu_2 = (-2.5, 2), \Sigma_1 = \Sigma_2 = \Sigma_3 = \text{diag}(0.5, 0.5)$ into 3D. Shown in Figure \ref{3Cyl_BLO_results}, the resulting obstacle looks roughly like three cylinders, and agent trajectories starting at different heights are identical when projected onto the first two dimensions.

\paragraph{Castle}

\begin{figure}[h]
\centering
\begin{minipage}{0.5\linewidth}
\includegraphics[width=.8\linewidth]{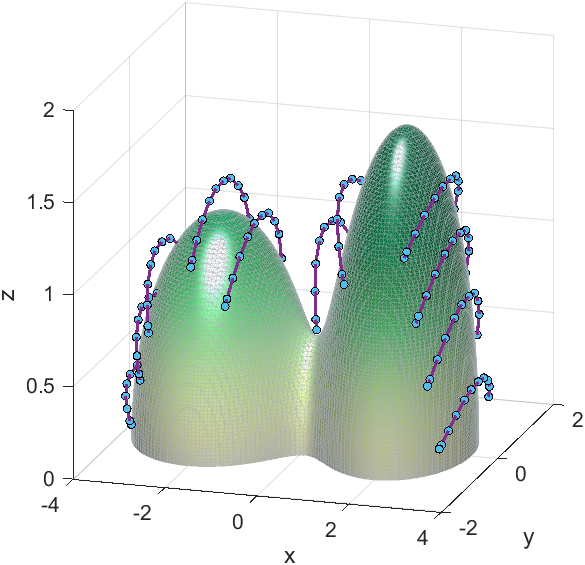}\\
\end{minipage}\hfill
\begin{minipage}{0.5\linewidth}
\includegraphics[width=.8\linewidth]{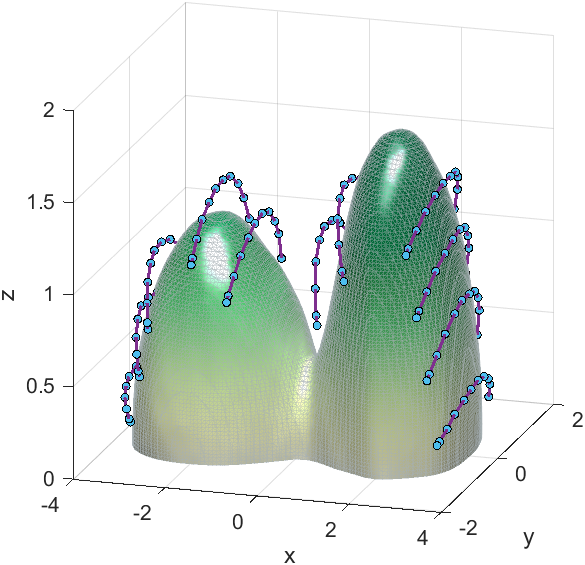}\\
\end{minipage}\hfill
\vskip -10pt
\centering
\begin{minipage}{0.5\linewidth}
\includegraphics[width=.8\linewidth]{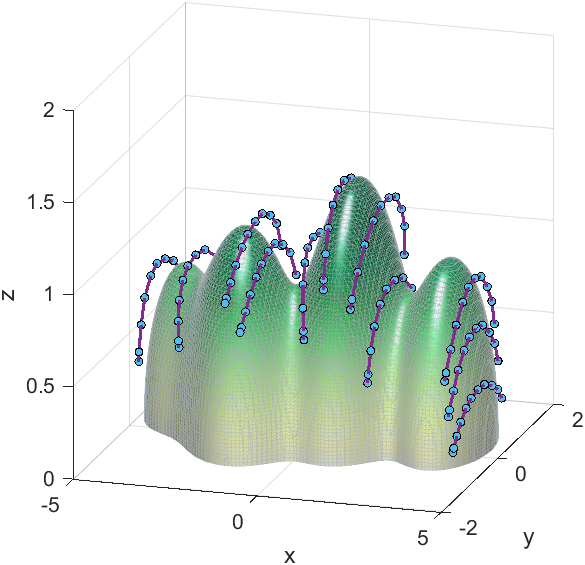}\\
\end{minipage}\hfill
\begin{minipage}{0.5\linewidth}
\includegraphics[width=.8\linewidth]{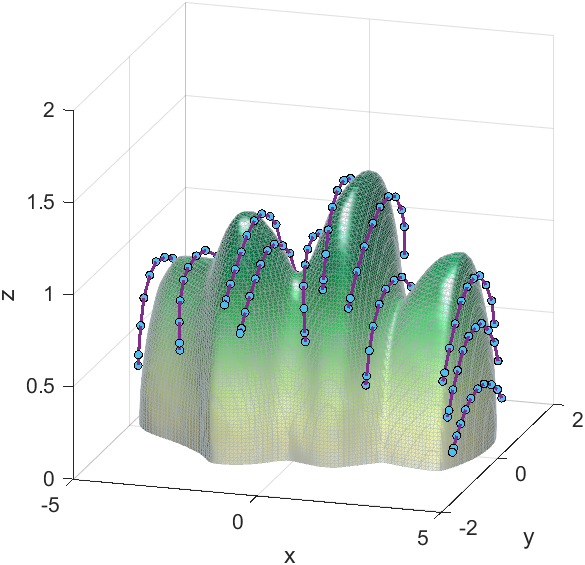}\\
\end{minipage}\hfill
\vskip -10pt
\caption{Left: ground truth castle(top) and mountain(bottom) obstacle with sample trajectories overlaid; Right: inferred obstacle(top) and mountain(bottom) from trajectories with predicted trajectories. The obstacle color gradient indicates its height. }
\label{Castle_BLO_results}
\end{figure}

Next, we consider obstacles with non-trivial patterns on the z-axis. In Figure \ref{Castle_BLO_results}, we present a mixture of two gaussians with different heights: $B(\vx) = 50\cdot \frac{1}{2}\sum_{i=1}^2 \mathrm{N}(\vx; \mu_i, \Sigma_i), \mu_1 = (-2, 0, 0), \mu_2 = (2, 0, 0), \Sigma_1 = \text{diag}(1.8, 0.8, 0.8), \Sigma_2 = \text{diag}(1, 0.6, 1.5)$. 

\paragraph{Mountains}

Lastly, we investigate a more complex obstacle resembling a mountain range in Figure \ref{Castle_BLO_results}. It is a mixture of four gaussians: $B(\vx) = 50\cdot \frac{1}{4}\sum_{i=1}^4 \mathrm{N}(\vx; \mu_i, \Sigma_i), \mu_1 = (1, 0.25, 0), \mu_2 = (3.5, 0.25, 0), \mu_3 = (-1.5, -0.35, 0), \mu_4 = (-4, 0.5, 0), \Sigma_1 = \text{diag}(0.5, 0.4, 0.8), \Sigma_2 = \text{diag}(1, 0.6, 0.8), \Sigma_3 = \text{diag}(1.2, 0.4, 1), \Sigma_4 = \text{diag}(1.5, 0.6, 1.2)$. \\

In each setup, we present the reconstructed obstacle and trajectories obtained with $N=30k$ training data and compare them to their ground truth counterparts. Since the obstacles are defined in all of $\mathbb{R}^3$, we visualize their iso-surfaces at appropriate values for intuitive inspection. Overall, our framework closely captures the obstacle iso-surfaces in all three scenarios. In addition, the learned trajectories qualitatively reflect agents' avoidance behaviors on the inferred obstacle.



\subsection{High Dimensional Gaussian}

\begin{figure}[h]
\centering
\begin{minipage}{0.25\linewidth}
\includegraphics[width=1\linewidth]{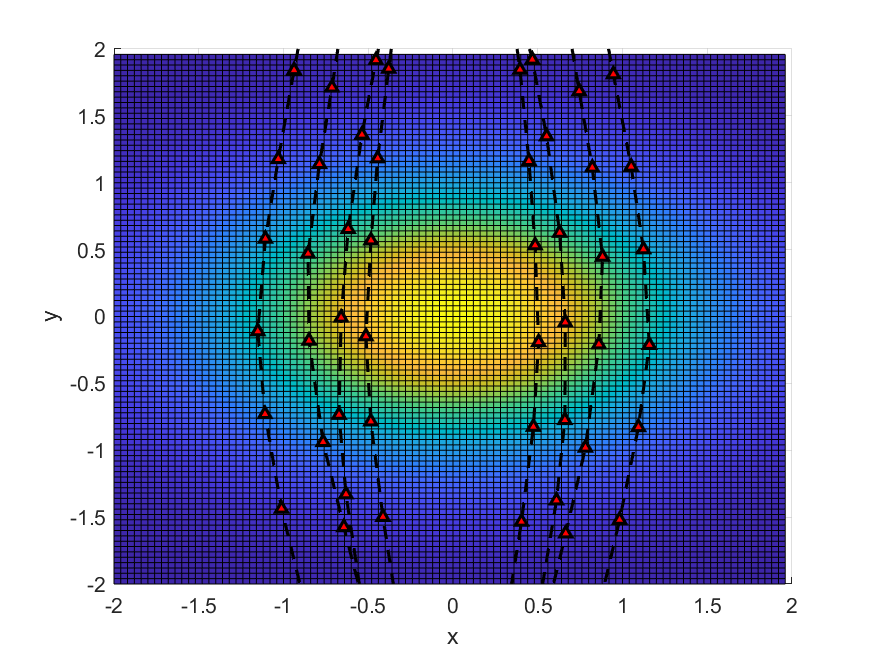}\\
\end{minipage}\hfill
\begin{minipage}{0.25\linewidth}
\includegraphics[width=1\linewidth]{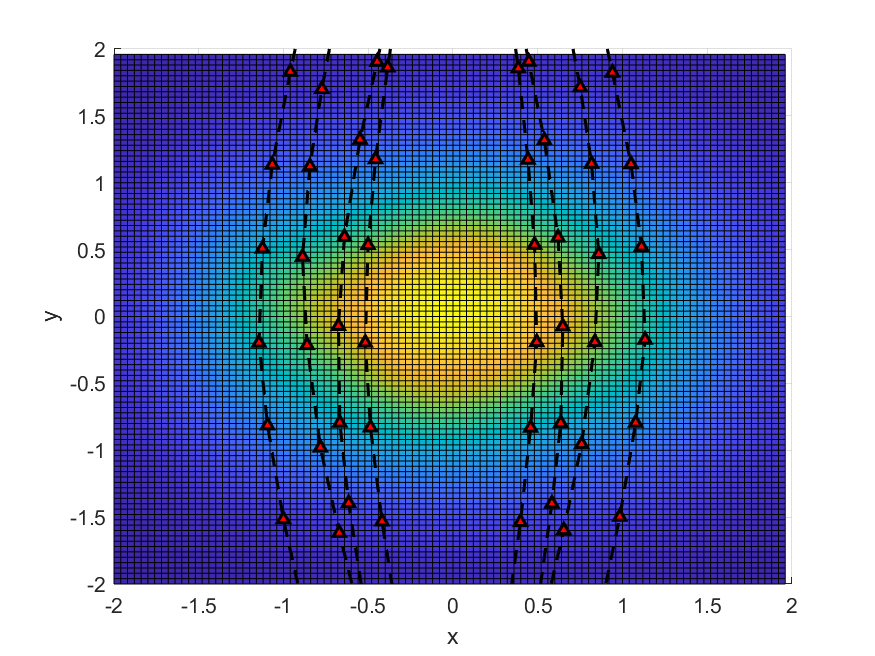}\\
\end{minipage}\hfill
\begin{minipage}{0.25\linewidth}
\includegraphics[width=1\linewidth]{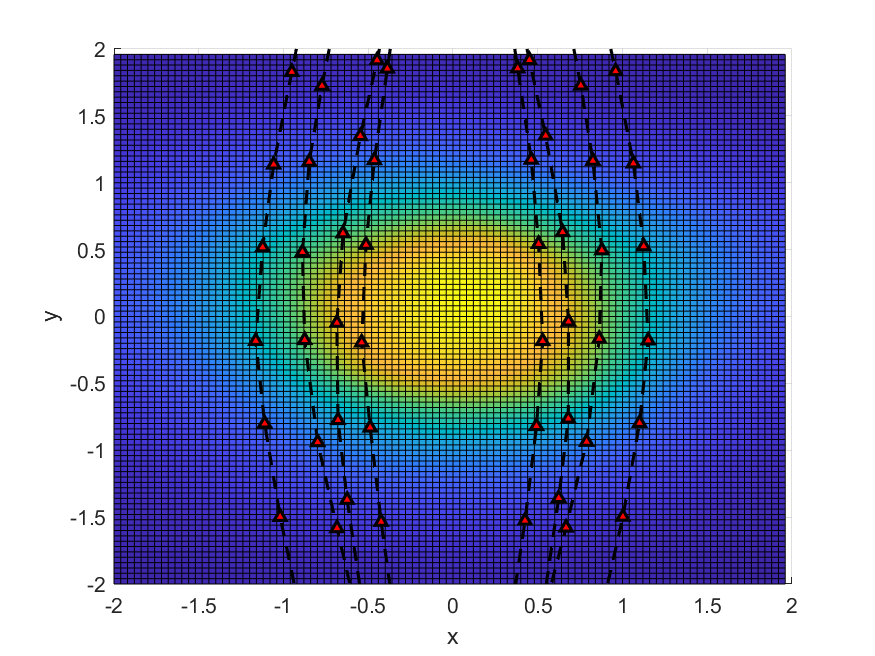}\\
\end{minipage}\hfill
\begin{minipage}{0.25\linewidth}
\includegraphics[width=1\linewidth]{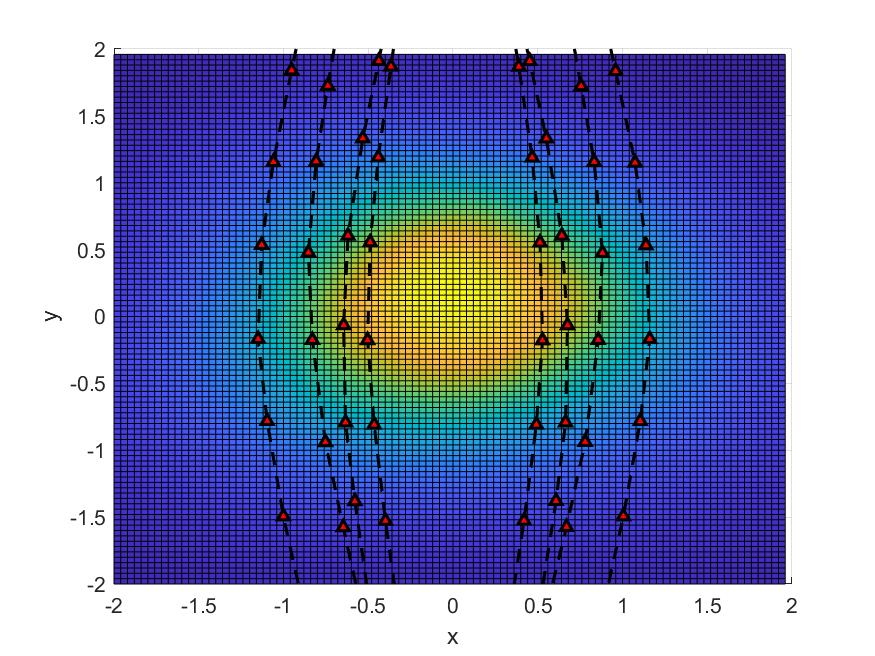}\\
\end{minipage}\hfill
\vskip -15pt
\caption{Left to right: Ground truth gaussian obstacle with sample trajectories overlaid; inferred obstacle from trajectories, with predicted sample trajectories overlaid in $d=2,5,10$ and projected onto the first two dimensions.}
\label{Gaussian_BLO_diff_dim_results}
\end{figure}

\begin{table}[h]
\begin{center}
\begin{small}
\begin{sc}
\begin{tabular}{lccccccc}
\toprule
Dataset & Gaussian, $d=2$ & Gaussian, $d=5$ & Gaussian, $d=10$\\
\midrule
Rel. $L_2$ err  & 0.071 & 0.062 & 0.088\\
\bottomrule
\end{tabular}
\end{sc}
\end{small}
\end{center}
\caption{Relative $L_2$ errors of recovered obstacles for different MFG settings.}
\label{table:obs_rel_l2_err_highd}
\vskip -15pt
\end{table}

A key advantage of employing NF parametrization for agent trajectories is the elimination of spatial discretization that is typically required for classic optimization methods \cite{achdou2010mean, benamou2014augmented,benamou2017variational,benamou2000computational, jacobs2019solving,papadakis2014optimal,yu2021fast}. Hence,  traditional methods are intractable in dimensions $d>3$ as the number of regular grid points scales exponentially. In contrast, our method sidesteps the curse of dimensionality and has the promise of inferring trajectories and obstacles in high dimensional spaces.

For empirical verification, we borrow a setting from \cite{NN_MFP} to construct a two-dimensional Gaussian obstacle situated between the initial and terminal distributions. The agents travel vertically with $p_0(\vx) = \mathrm{N}(\vx;(0, 3), 0.3I), p_1(\vx) = \mathrm{N}(\vx; (0, -3), 0.3I)$, and the obstacle is $B(\vx) = 50 \cdot \mathrm{N} (\vx; 0, \text{diag}(1,0.5))$. For $d>2$, we extend the gaussian constantly beyond the first two dimensions and seek to recover the obstacle in $d=5,10$ dimensions. 

In our setup, the optimal agent trajectories and true obstacles are identical when projected to the first two dimensions regardless of the problem dimension. This behavior is captured accurately in  figure~\ref{Gaussian_BLO_diff_dim_results}, so our framework faithfully recovers trajectories and obstacles across various dimensions. In addition, our method has similar runtime speeds across different dimensions thanks to GPU parallelism.


\begin{figure}[th]
\vskip -5pt
\centering
\begin{minipage}{0.25\linewidth}
\includegraphics[width=1\linewidth]{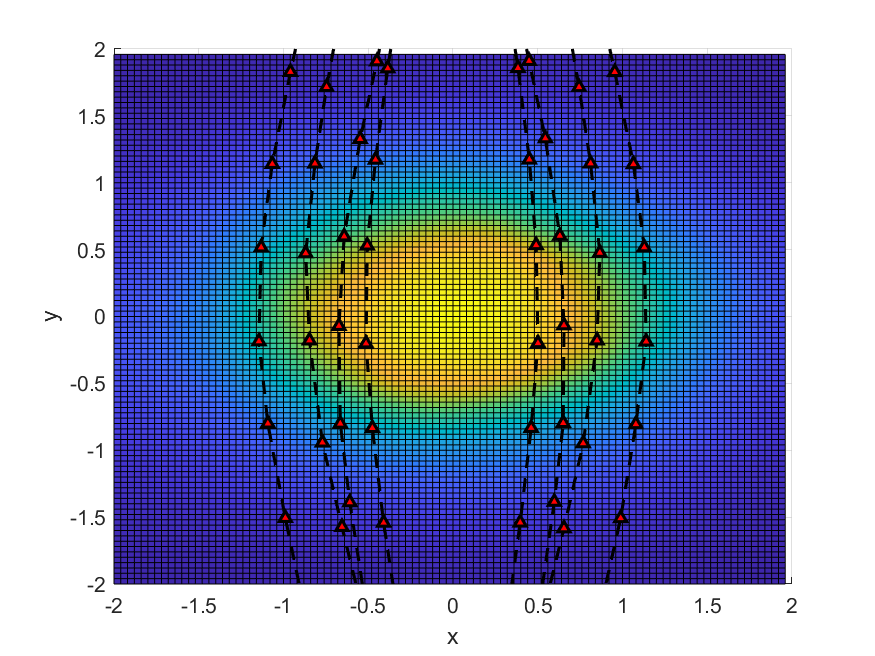}\\
\includegraphics[width=1\linewidth]{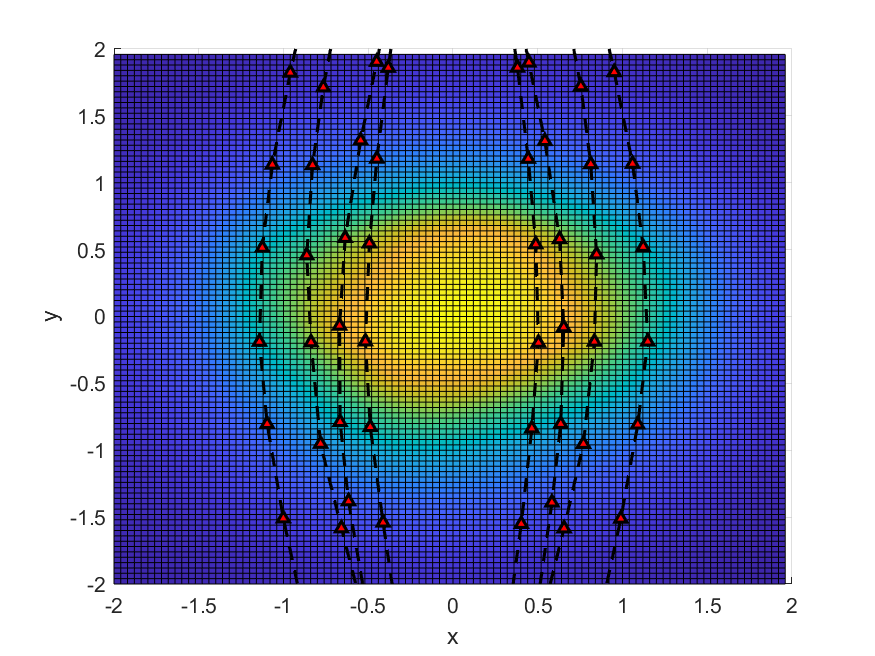}\\
\end{minipage}\hfill
\begin{minipage}{0.25\linewidth}
\includegraphics[width=1\linewidth]{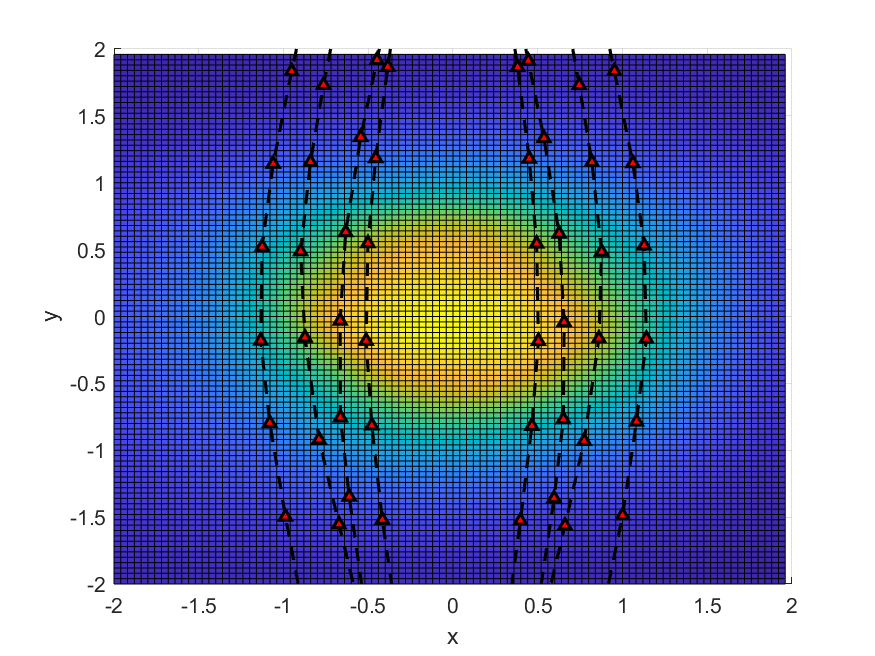}\\
\includegraphics[width=1\linewidth]{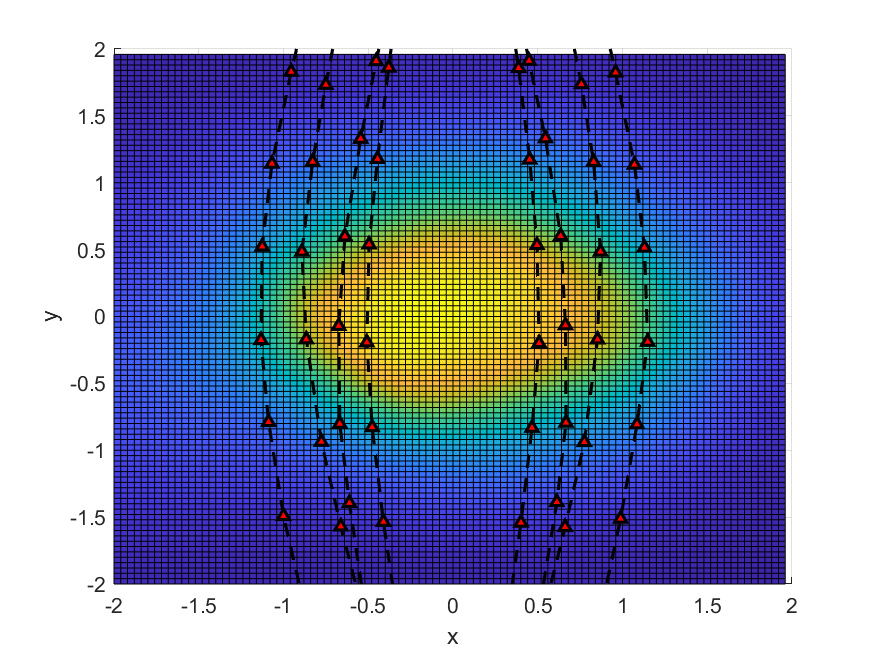}\\
\end{minipage}\hfill
\begin{minipage}{0.5\linewidth}
\includegraphics[width=1\linewidth]{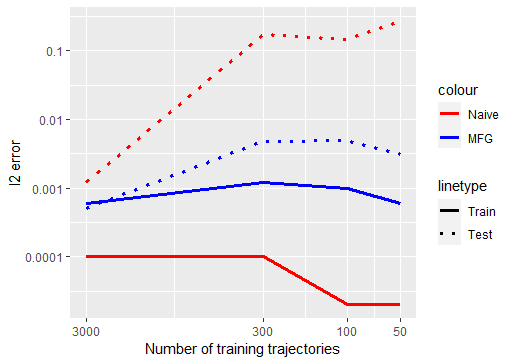}\\
\end{minipage}\hfill
\vskip -10pt
\caption{Left, counter-clockwise: inferred obstacles and trajectories from $N=3000, 300, 100, 50$ training data. Right: comparing training and testing performance for the naive (minimizing $l_2$ trajectory error) and MFG approaches as data. It is evident that the MFG approach generalizes well even when trained on very little data, while the naive method overfits more as data become more scarce.}
\label{fig:Gaussian_BLO_low_data}
\end{figure}

\subsection{Learning with Scarce Data}

\begin{figure}[H]
\centering
\begin{minipage}{0.25\linewidth}
\includegraphics[width=1\linewidth]{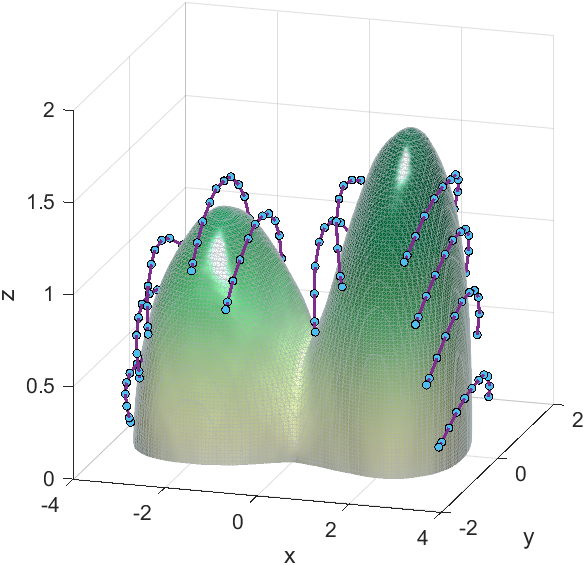}\\
\includegraphics[width=1\linewidth]{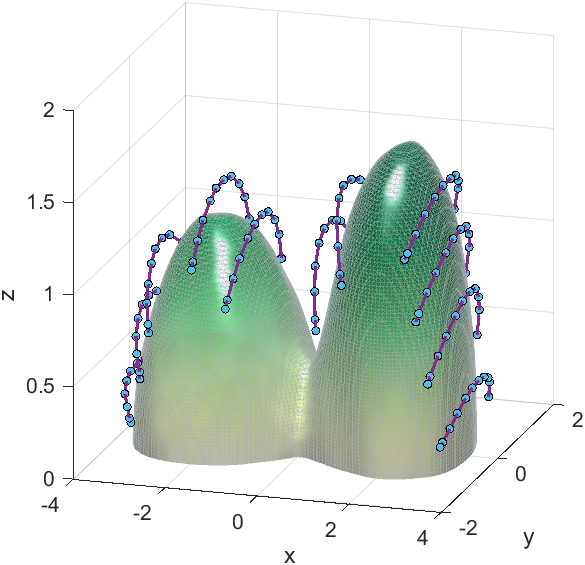}\\
\end{minipage}\hfill
\begin{minipage}{0.25\linewidth}
\includegraphics[width=1\linewidth]{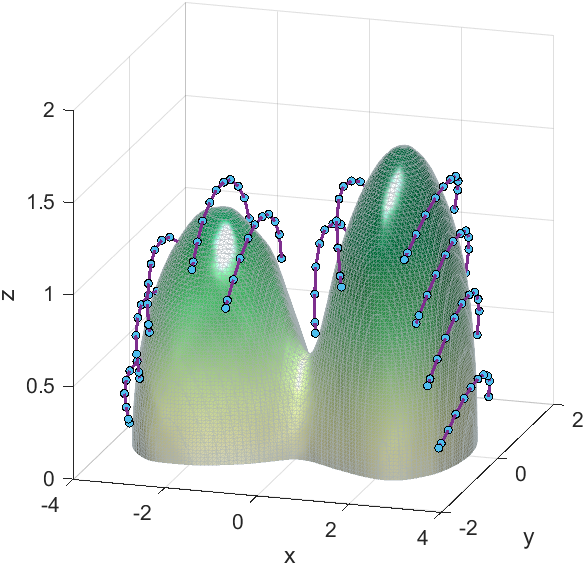}\\
\includegraphics[width=1\linewidth]{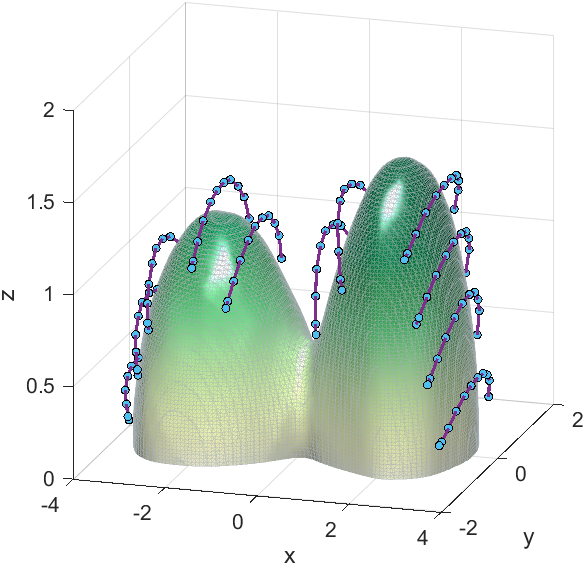}\\
\end{minipage}\hfill
\begin{minipage}{0.5\linewidth}
\includegraphics[width=1\linewidth]{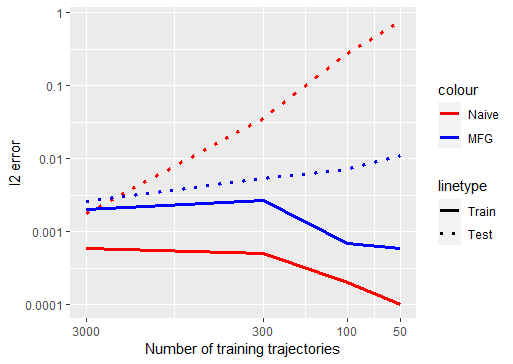}\\
\end{minipage}\hfill
\vskip -10pt
\caption{Left, counter-clockwise: inferred obstacles and trajectories from $N=3000, 300, 100, 50$ training data. Right: comparing training and testing performance for the naive (minimizing $l_2$ trajectory error) and MFG approaches as data. It is evident that the MFG approach generalizes well even when trained on very little data, while the naive method overfits more as data become more scarce.}
\label{fig:castle_BLO_low_data}
\end{figure}

In reality, obtaining trajectory data may be costly in both time and resources. As a result, it is highly desirable if our approach wagers a reasonable guess at the latent obstacle even if few training trajectories are available. In the following experiments, we explore the possibility of obstacle recovery in the low data regime for different settings. Surprisingly, our approach performs reasonably well even when the size of the training set is reduced by more than two orders of magnitude: $N \sim 10^4 \to N \sim 10^2$, as evidenced in left curves in Figures \ref{fig:Gaussian_BLO_low_data} and \ref{fig:castle_BLO_low_data}.

Furthermore, our framework demonstrates great efficiency for trajectory learning even if obstacle detection is of secondary interest. While a simple approach of fitting a network on the $l_2$ error of trajectories may suffice when there is an abundance of data, it becomes less effective when data is scarce, and additional assumptions about the underlying dynamics are needed. By balancing a naive likelihood modeling on sampled trajectories with a prior that reflects their mean-field nature, the cost-minimization constraint provides a strong inductive bias that constrains the hypothesis space of agent trajectories, thereby enabling tractable learning with very few training data. 

To corroborate this intuition empirically, we plot and compare the trajectory generalization errors for different number of training data between the naive approach and the inverse MFG framework. Figure~\ref{fig:Gaussian_BLO_low_data} shows that the naive approach fits the training trajectories well but struggles to generalize when data is limited. In contrast, our approach reflects the true dynamics much more accurately in the scarce data regime. In other words, the inverse MFG framework serves as an implicit regularization to effectively reduce overfitting on the learned trajectories. 

In addition, we replicate the above settings on the castle obstacle, which a distribution with many modes supported in 3D. The observed trends shown in figure~\ref{fig:castle_BLO_low_data} are qualitatively identical, which speaks to the robustness of our approach.

\section{Conclusion}

In this work, we study an inverse MFG problem to simultaneously infer latent obstacles and optimal agent trajectories given partial observation on the latter. To this end, we formulate and solve a bilevel optimization problem with the penalty method in a variety of MFG setups. Our method exhibits a convincing degree of robustness across different obstacle complexity and dimensionality in synthetic setups, which suggests potential future applications in e.g., obstacle recovery in configuration spaces. In addition, our method serves as a MFG regularizor on maximum likelihood trajectory learning, which improves generalization particularly when the available training data is limited.

\section{Acknowledgment}
R. Lai and T. Chen's research is supported in part by NSF SCALE MoDL (DMS-2401297). 

\bibliography{references}
\bibliographystyle{plain}



\end{document}